\documentclass{amsart}
\usepackage[mathscr]{eucal}
\usepackage{graphicx}
\usepackage{amscd}
\usepackage{amsmath}
\usepackage{amsthm}
\usepackage{amsxtra}
\usepackage{calc}
\usepackage{amsfonts}
\usepackage{amssymb}
\usepackage{latexsym}
\usepackage{pdfsync}
\usepackage{amsopn}
\usepackage{mathrsfs}
\usepackage{accents}
\usepackage[all]{xy}
\SelectTips{cm}{}
\calclayout

\usepackage[colorlinks=true,linkcolor=red,citecolor=green,bookmarks=true]{hyperref}

\newcommand{\new}[1]{{\color{blue} #1}}

\newcommand{\rt}{\rightarrow}
\newcommand{\lra}{{\xymatrix@C=3em{\ar @{>} [r] & }}}
\newcommand{\lrt}{\longrightarrow}

\def\db{\operatorname{\mathsf{D^b}}}

\def\ds{\operatorname{\mathsf{D_{sg}}}}

\newcommand{\md}{\mathsf{mod}}
\newcommand{\si}{\mathsf{\Sigma}}
\newcommand{\syz}{\mathsf{\Omega}}

\newcommand{\st}{\stackrel}

\newcommand{\HT}{\mathsf{H}}

\newcommand{\Z}{\mathbb{Z} }

\newcommand{\uEnd}{\underline{\mathsf{End}}}
\newcommand{\ugp}{\underline{\mathsf{Gp}}}
\newcommand{\gp}{{\mathsf{Gp}}}
\newcommand{\n}{{\mathfrak{n}}}

\newcommand{\p}{{\mathfrak{p}}}
\newcommand{\cok}{{\rm{Coker}}}

\newcommand{\Tr}{\mathrm{Tr}}

\newcommand{\di}{\mathsf{dim}}

\newcommand{\Tor}{\mathsf{Tor}}

\newcommand{\Hom}{{\mathsf{Hom}}}

\newcommand{\End}{{\mathsf{End}}}

\newcommand{\id}{{\mathsf{id}}}

\newcommand{\mon}{{\mathsf{Mon}}}
\newcommand{\mor}{{\mathsf{Mor}}}

\def\si{\operatorname{\mathsf{\Sigma}}}

\newcommand{\Ext}{\mathsf{{Ext}}}

\newcommand{\A}{\mathcal{{A}}}

\newcommand{\fo}{\footnotesize }

\newcommand{\cp} {\mathcal{P}}

\newtheorem{theorem}{Theorem}[section]

\newtheorem{cor}[theorem]{Corollary}

\newtheorem{lemma}[theorem]{Lemma}
\newtheorem{prop}[theorem]{Proposition}

\theoremstyle{definition}
\newtheorem{example}[theorem]{Example}

\newtheorem{remark}[theorem]{Remark}

\newtheorem{s}[theorem]{}

\theoremstyle{plain}

\theoremstyle{definition}
\newtheorem{dfn}[theorem]{Definition}

\numberwithin{equation}{section}

\begin{document}

\title[The homotopy category of monomorphisms between projective modules]
{The homotopy category of monomorphisms between projective modules}

\author[Bahlekeh, Fotouhi, Nateghi and Salarian]{Abdolnaser Bahlekeh, Fahimeh Sadat Fotouhi, Armin Nateghi and Shokrollah Salarian}

\address{Department of Mathematics, Gonbad Kavous University, Postal Code:4971799151, Gonbad Kavous, Iran}
\email{bahlekeh@gonbad.ac.ir}
\address{School of Mathematics, Institute for Research in Fundamental Science (IPM), P.O.Box: 19395-5746, Tehran, Iran}
\email{ffotouhi@ipm.ir}

\address{Department of Mathematics, University of Isfahan, P.O.Box: 81746-73441, Isfahan,
Iran }\email{a.nateghi92@sci.ui.ac.ir}

\address{Department of Mathematics, University of Isfahan, P.O.Box: 81746-73441, Isfahan,
Iran and \\ School of Mathematics, Institute for Research in Fundamental Science (IPM), P.O.Box: 19395-5746, Tehran, Iran}
\email{Salarian@ipm.ir}

%\date{\today, \setcounter{hours}{\time/60} \setcounter{minutes}{\time-\value{hours}*60} \thehours\,h\ \theminutes\,min}

\subjclass[2020]{13D09, 18G80, 16G70, 16G30}

\keywords{monomorphism category, homotopy category, almost split sequence, Auslander-Reiten translation, singularity category.}

\thanks{{This work is based upon research funded by Iran National Science Foundation (INSF) under project No. 4001480.} The research of the second author was in part supported by a grant from IPM}

\begin{abstract}
Let $(S, \n)$ be a commutative noetherian local ring and $\omega\in\n$ be non-zerodivisor. This paper deals with the behavior of the category $\mon(\omega, \cp)$ consisting of all monomorphisms between finitely generated projective $S$-modules with cokernels annihilated by $\omega$. We introduce a homotopy category $\HT\mon(\omega, \cp)$, which is shown to be triangulated. It is proved that this homotopy category embeds into the singularity category of the factor ring  $R=S/{(\omega)}$. As an application, not only the existence of almost split sequences {ending at indecomposable non-projective objects of} $\mon(\omega, \cp)$ is proven, but also the Auslander-Reiten translation, $\tau_{\mon}(-)$, is completely recognized. Particularly, it will be observed that any non-projective object of $\mon(\omega, \cp)$ with local endomorphism ring is invariant under the square of the Auslander-Reiten translation.
\end{abstract}

\maketitle

%\tableofcontents
\section{Introduction}

Let $\Lambda$ be an associative ring with identity and $\md\Lambda$ the category of finitely generated left $\Lambda$-modules. The morphism category $\mor(\Lambda)$ of $\Lambda$ has as objects the maps in $\md\Lambda$ and whose morphisms are given by commutative squares. The monomorphism category $\mon(\Lambda)$ of $\Lambda$ is the full subcategory of $\mor(\Lambda)$ consisting of all monomorphisms, which is also known as the submodule category. It is known that $\mor(\Lambda)$ is an abelian category (\cite[Proposition 1.1]{fgr}) and $\mon(\Lambda)$ is an extension-closed additive subcategory of $\mor(\Lambda)$, and so, it will be an exact category in the sense of Quillen \cite[Appendix A]{kel}.

The study of monomorphism categories goes back to Birkhoff \cite{bir} in 1934, who initiated classifying the indecomposable objects of the category $\mon(\Z/{(p^t)})$ with $t\geq 2$ and $p$ a prime number. It is well-understood that the categories $\mor(\Lambda)$ and $\mon(\Lambda)$ are much more complicated than underlying module categories $\md\Lambda$. This subject has been recently attracted more attention and
studied in a systematic and deep work by Ringel and Schmidmeier \cite{RS1, RS2}. See also \cite{ha21, hz, luo2017separated, xiong2014auslander, zhang2019separated} for a generalization of such works. Particularly, Kussin, Lenzing and Meltzer \cite{kussin2013nilpotent}
have discovered a surprising link between the stable submodule category with the singularity theory via weighted projective lines of type (2, 3, $p$).

A nice result of Ringel and Schmidmeier
indicates that, over an artin algebra $\Lambda$, the Auslander-Reiten translation in the monomorphism category
can be computed within $\md\Lambda$ by using their construction of minimal monomorphisms. To be more precise, for a given object $(M\st{f}\rt N)\in\mor(\Lambda)$,
they define the minimal monomorphism $\mathsf{Mimo}(f)$ of $f$ as $\mathsf{Mimo}(f)=[f~~e]:M\rt N\oplus E$, where $e$ is an extenstion of an injective envelope $e':\ker(f)\rt E$. In a remarkable result, they
have shown that, if $f$ is an indecomposable object of $\mon(\Lambda)$, then $\tau_{\mon}(f)=\mathsf{Mimo}\tau_{\Lambda}\cok(f)$;  see \cite[Theorem 5.1]{RS2}. Here $\tau_{\mon}(f)$ is the Auslander-Reiten
translation of $f$ in $\mon(\Lambda)$. This result has been generalized over noetherian algebras in \cite[Theorem 5.9]{bfs}.

From now on, assume that $(S, \n)$ is a commutative noetherian local ring with $\dim S\geq 2$ and $\omega\in\n$ is non-zerodivisor. Assume that $\mon(\omega, \cp)$ is the full subcategory of $\mon(S)$ consisting of all monomorphisms $(P\st{f}\rt Q)$ in the module category $\md S$ such that $P$ and $Q$ are finitely generated projective modules and $\cok f$ is annihilated by $\omega$. In this paper, we will show that this category is well-behaved. Firstly, we need to give some definitions. A morphism  $\psi=(\psi_1, \psi_0):(P\st{f}\rt Q)\lrt (P'\st{f'}\rt Q')$ in $\mon(\omega, \cp)$ is called  null-homotopic, if there are  $S$-homomorphisms $s_1:P\rt Q'$ and $s_0:Q\rt P'$ such that $f'\psi_1-f's_0f=\omega. s_1$, or equivalently, $\psi_0f-f's_0f=\omega. s_1$.  Now we define the homotopy category $\HT\mon(\omega, \cp)$ with the same objects as $\mon(\omega, \cp)$ and its morphism sets are morphism sets in $\mon(\omega, \cp)$ modulo null-homotopic.  It is fairly easy to see that for a given object $(P\st{f}\rt Q)\in\mon(\omega, \cp)$, there is a unique morphism $Q\st{f_{\si}}\rt P$ such that $ff_{\si}=\omega.\id_Q$ and $f_{\si}f=\omega.\id_P$, and in particular, $(Q\st{f_{\si}}\rt P)$ is also an object of $\mon(\omega, \cp)$; see Lemma \ref{lem1}. It is proved that $\si:\HT\mon(\omega, \cp)\lrt\HT\mon(\omega, \cp)$, with $\si(f)=-f_{\si}$  and $\si(\psi_1, \psi_0)=(\psi_0, \psi_1)$, is an auto-equivalence functor, and particularly, with $\si$ being suspension, $\HT\mon(\omega, \cp)$ gets a triangulated structure in a natural way; see Proposition \ref{proptri}. Also a tie connection between the homotopy category $\HT\mon(\omega, \cp)$ and the singularity category of the factor ring $R=S/{(\omega)}$, will be discovered. Precisely, it will be shown that there is a fully faithful triangle functor $\HT\mon(\omega, \cp)\lrt\ds(R)$; see Corollary \ref{ccc}. By the aid of this result, in the last section, we will able to show that each indecomposable non-projective object of the category $\mon(\omega, \cp)$ appears as the right term of an almost split sequence, provided that $R$ is a complete Gorenstein ring which is an isolated singularity. Particularly, an explicit description of the Auslander-Reiten translation in {$\mon(\omega, \cp)$, which will be also denoted by $\tau_{\mon}(-)$,} is given. Precisely, it is shown that for a given non-projective indecomposable object $(P\st{f}\rt Q)\in\mon(\omega, \cp)$, we have $\tau_{\mon}(f)=f$, if $\di R$ is even, and otherwise, $\tau_{\mon}(f)=f_{\si}$; see Theorem \ref{ass}. This particularly implies that $\tau^2_{\mon}(f)=f$. This result should be compared with Corollary 6.5 of \cite{RS2}, where they have shown that if $\Lambda$ is a commutative uniserial algebra, then for an indecomposable non-projective object $(M\st{f}\rt N)$ in $\mon(\Lambda)$, there is an isomorphism of
objects $\tau_{\mon}^6(f)\cong f$.

Throughout the paper, $(S, \mathfrak{n})$ is a commutative noetherian local ring with maximal ideal $\n$, $\omega\in\n$ is a non-zerodivisor element, and $R$ is the factor ring $S/{(\omega)}$. Unless otherwise specified, by a module we mean a finitely generated $S$-module and $\md S$ stands for the category of all finitely generated $S$-modules. Moreover, to consider a map $f:M\rt N$ as an object of $\mor(S)$ we use parentheses and denote it by $(M\st{f}\rt N)$.

\section{The homotopy category of monomorphisms}
This section is devoted to introduce and study the homotopy category of the monomorphism category of projective $S$-modules. Among others, we show that this subcategory admits a triangulated structure. Moreover, it is proved that it can be considered as a full triangulated subcategory of the singularity category $\ds(R)$ of $R$. We begin with the following definition.

\begin{dfn}
By the category $\mon(\omega, \cp)$, we mean a category that whose objects are those $S$-monomorphisms $(P\st{f}\rt Q)$, where $P, Q\in\cp(S)$ and $\cok f$ is an $R$-module. Here $\cp(S)$ is the category of all (finitely generated) projective $S$-modules. Moreover, a morphism $\psi=(\psi_1, \psi_0):(P\st{f}\rt Q)\lrt (P'\st{f'}\rt Q')$ between two objects is a pair of $S$-homomorphisms $\psi_1:P\rt P'$ and $\psi_0:Q\rt Q'$ such that $\psi_0f=f'\psi_1$.  It is clear that $\mon(\omega, \cp)$ is a full additive subcategory of the monomorphism category $\mon(S)$.
% the square $$\xymatrix{P\ar[r]^{f} \ar[d]_{\psi_1} &\ Q \ar[d]_{\psi_0}\\ P' \ar[r]^{f'} & Q'}$$commutes. It is clear that $\mon(\omega, \cp)$ is a full additive subcategory $\mon(S)$.
\end{dfn}

\begin{dfn}
We say that a morphism  $\psi=(\psi_1, \psi_0):(P\st{f}\rt Q)\lrt (P'\st{f'}\rt Q')$ in $\mon(\omega, \cp)$ is {\em null-homotopic}, if there are  $S$-homomorphisms $s_1:P\rt Q'$ and $s_0:Q\rt P'$ such that $f'\psi_1-f's_0f=\omega. s_1$, or equivalently, $\psi_0f-f's_0f=\omega. s_1$.
\end{dfn}

The {\em homotopy category} $\HT\mon(\omega, \cp)$ of $\mon(\omega, \cp)$, is defined as follows; its objects are the same as $\mon(\omega, \cp)$ and its morphism sets are morphism sets in $\mon(\omega, \cp)$ modulo null-homotopic. It is easily seen that null-homotopies are compatible with addition and composition of morphisms in $\mon(\omega, \cp)$. This, in conjunction with the fact that $\mon(\omega, \cp)$ is an additive category, would imply the result below.

\begin{prop}The homotopy category  $\HT\mon(\omega, \cp)$  is an additive category.
\end{prop}

We should emphasize that although the proof of the next result is the same as for Lemma 2.6  of \cite{bfns}, we include it only for the sake of completeness.

\begin{lemma}\label{lem1}
Let $(P\st{f}\rt Q)\in\mon(\omega, \cp)$ be arbitrary. Then there exists
a unique morphism $(Q\st{f_{\si}}\rt P)$ such that $ff_{\si}=\omega.\id_{Q}$ and $f_{\si}f=\omega.\id_{P}$. In particular, $(Q\st{f_{\si}}\rt P)$ is an object of $\mon(\omega, \cp)$.
\end{lemma}
\begin{proof}By our assumption, there is a short exact sequence of $S$-modules; $0\rt P\st{f}\rt Q\rt\cok f\rt 0$ such that {$\omega\cok f=0$. Take the following commutative diagram with exact rows;
\[\xymatrix{&\\ \eta:0 \ar[r] & P \ar[r]^{f} \ar[d]_{\omega} & Q \ar[r] \ar[d]_{\omega} & \cok f \ar[r] \ar[d]_{\omega} & 0
\\ \eta: 0 \ar[r] & P \ar[r]^{f} & Q \ar[r] & \cok f \ar[r] & 0.}\]So, applying \cite[Lemma 1.1, page 163]{mit} gives us the following commutative diagram with exact rows;
\[\xymatrix{&\\ \eta:0 \ar[r] & P \ar[r]^{f} \ar[d]_{\omega} & Q \ar[r] \ar[d] & \cok f \ar[r] \ar@{=}[d] & 0
\\  0 \ar[r] & P\ar@{=}[d] \ar[r] & T \ar[r]\ar[d] & \cok f \ar[r]\ar[d]_{\omega} & 0\\  \eta:0 \ar[r] & P \ar[r]^{f} & Q \ar[r] & \cok f \ar[r] & 0.}\]Since $\omega\cok f=0$, the middle row will be split,} and so, there is a morphism $f_{\si}:Q\rt P$ such that $f_{\si}f=\omega.\id_{P}$. Another use of the fact that $\omega\cok f=0$, leads us to infer that $\omega Q\subseteq f(P)$. This fact besides the equality $f_{\si}f=\omega.\id_{P}$ would imply that $ff_{\si}=\omega.\id_{Q}$. It should be noted that, if there is another morphism $g:Q\rt P$ satisfying the mentioned conditions, then we will have $ff_{\si}=fg$, and then, $f$ being a monomorphism ensures the validity of the equality $f_{\si}=g$.
Now we show that $(Q\st{f_{\si}}\rt P)\in\mon(\omega, \cp)$. As $f_{\si}$ is evidently a monomorphism, we only need to check that $\omega$ annihilates $\cok f_{\si}$. To see this, consider the short exact sequence $0\rt Q\st{f_{\si}}\rt P\st{\pi}\rt\cok f_{\si}\rt 0$. For a given object $y\in\cok f_{\si}$, take $x\in P$ such that $\pi(x)=y$. Since $f_{\si}f=\omega.\id_{P}$, we have $\omega P\subseteq f_{\si}(Q)$.  Consequently, $\omega y=\omega\pi(x)=\pi(\omega x)\in\pi f_{\si}(Q)=0$, meaning that $\omega\cok f_{\si}=0$, as needed.
\end{proof}

\begin{cor}For a given object $(P\st{f}\rt Q)\in\mon(\omega, \cp)$, $(f_{\si})_{\si}=f$. 
\end{cor}

\begin{remark}\label{inv}Assume that $\psi=(\psi_1, \psi_0):(P\st{f}\rt Q)\lrt(P'\st{f'}\rt Q')$ is a morphism in $\mon(\omega, \cp)$.  Since $f'\psi_1=\psi_0f$, we will  have the equalities;   $f'\psi_1f_{\si}=\psi_0ff_{\si}=\omega.\psi_0=f'f'_{\si}\psi_0$. Now  $f'$ being a monomorphism, gives rise to the equality $\psi_1f_{\si}=f'_{\si}\psi_0$. Namely $(\psi_0, \psi_1):(Q\st{f_{\si}}\rt P)\lrt (Q'\st{f'_{\si}}\rt P')$ is also a morphism in the category $\mon(\omega, \cp)$. Particularly, in a similar way, one may see that the morphisms $(\psi_1, \psi_0)$ and $(\psi_0, \psi_1)$ are null-homotopic, simultaneously.
\end{remark}

In what follows, we intend to show that the category $\HT\mon(\omega, \cp)$, admits a natural structure of triangulated category. In this direction,  we need to determine a translation functor $\si$ and a class of exact triangles.

Assume that $(P\st{f}\rt Q)$ is an arbitrary object of $\mon(\omega, \cp)$. In view of Lemma \ref{lem1}, there is a unique object $(Q\st{f_{\si}}\rt P)\in\mon(\omega, \cp)$. Now we define $\si((P\st{f}\rt Q)):=(Q\st{-f_{\si}}\rt P)$. Moreover, for a given morphism $\psi=(\psi_1, \psi_0)$ in  $\HT\mon(\omega, \cp)$, we set $\si((\psi_1, \psi_0)):=(\psi_0, \psi_1)$.
So, one may easily see that $\si:\HT\mon(\omega, \cp)\lrt\HT\mon(\omega, \cp)$ is an additive functor, with $\si^{2}$ the identity functor, because by Lemma \ref{lem1}, $\si(P\st{f}\rt Q)$ is unique.  Precisely, we have the result below.
\begin{prop}$\si:\HT\mon(\omega, \cp)\lrt\HT\mon(\omega, \cp)$ is an auto-equivalence functor.
\end{prop}

Assume that  $\psi=(\psi_1, \psi_0):(P\st{f}\rt Q)\lrt (P'\st{f'}\rt Q')$  is a morphism in  $\HT\mon(\omega, \cp)$.  According to Lemma \ref{lem1}, there are objects $(Q\st{f_{\si}}\rt P)$ and $(Q'\st{f'_{\si}}\rt P')$  in $\mon(\omega, \cp)$ such that $ff_{\si}=\omega.\id_{Q}$ and $f'f'_{\si}=\omega.\id_{Q'}$. Now we define the {mapping cone} $C(\psi)$ of $\psi$, as $C(\psi):=(P'\oplus Q\st{c}\rt Q'\oplus P)$, where $c=\tiny {\left[\begin{array}{ll} f' & \psi_0 \\ {0} & {-f_{\si}} \end{array} \right].}$  As we will observe below, $C(\psi)$ is an object of $\mon(\omega, \cp)$. Moreover, we have that  $c_{\si}=\tiny {\left[\begin{array}{ll} f'_{\si} & \psi_1 \\ {0} & {-f} \end{array} \right].}$

\begin{lemma}\label{tt}With the notation above, $C(\psi)\in\mon(\omega, \cp)$.
\end{lemma}
\begin{proof}
 Consider the exact sequence $0\lrt(P'\st{f'}\rt Q')\lrt(P'\oplus Q\st{c}\rt Q'\oplus P)\lrt(Q\st{-f_{\si}}\rt P)\lrt 0$ with split rows. Since $f'$ and $f_{\si}$ are monomorphisms, the same will be true for $c$. So, it remains to show that $\omega\cok c=0$. Take an object $(x, y)\in Q'\oplus P$. We intend to find an object $(u, z)\in P'\oplus Q$ with $c(u, z)=(\omega x, \omega y)$. As $\omega\cok f_{\si}=0$, one may find an object $z\in Q$ such that $-f_{\si}(z)=\omega y$. Similarly, there is an object $t\in P'$ in which $f'(t)=\omega x$. So, applying Remark \ref{inv} gives us the equality  $-f'_{\si}\psi_0(z)=\omega\psi_1(y)$. Thus, we will have the equality $-f'f'_{\si}\psi_0(z)=\omega f'\psi_1(y)$. Now, since $f'f'_{\si}=\omega.\id_{Q'}$ and $\omega$ is non-zerodivisor, one may deduce that $-\psi_0(z)=f'\psi_1(y)$. Hence, putting $u:=t+\psi_1(y)$, we will have $c(u, z)=(\omega x, \omega y)$. This, indeed, means that $\omega\cok c=0$, and so, the proof is finished.
\end{proof}

Now we define a standard triangle in the category $\HT\mon(\omega, \cp)$  as a triangle of the form, $(P\st{f}\rt Q)\st{\psi}\lrt (P'\st{f'}\rt Q')\st{{[\id~~0]^t}}\lrt C(\psi)\st{{[0~~\id]}}\lrt {\si(P\st{f}\rt Q)}{=(Q\st{-f_{\si}}\rt P)}$. We say that a triangle in  $\HT\mon(\omega, \cp)$ is an {\em exact triangle}, if it is isomorphic to a standard triangle. Assume that $\Delta$ is the collection of all exact triangles in the homotopy category  $\HT\mon(\omega, \cp)$.{ Before proving that $\HT\mon(\omega, \cp)$ admits a triangulated structure, we need to state some preliminary results.

The proof of the next result follows by the same argument given in \cite[Lemma 2.14]{bfns}, and we include its proof for the sake of completeness.
\begin{lemma}\label{proj}Let $Q$ be a projective $S$-module. Then $(Q\st{\id}\rt Q)$ and $(Q\st{\omega}\rt Q)$ are projective objects in $\mon(\omega, \cp)$.
\end{lemma}
\begin{proof}We deal only with the case  $(Q\st{\omega}\rt Q)$, because the other one is obtained easily. Take a short exact sequence $0\lrt (E_1\st{e_1}\rt E_0)\lrt (T_1\st{g_1}\rt T_0)\st{\varphi}\lrt (Q\st{\omega}\rt Q)\lrt 0$ in $\mon(\omega, \cp)$. Now projectivity of  $Q$  gives us a morphism $\psi_0:Q\rt T_0$ with $\varphi_0\psi_0=\id_Q$. Since $\cok g_1$ is annihilated by $\omega$, one may find a morphism $\psi_1:Q\rt T_1$ making the following diagram commutative
{\footnotesize\[\xymatrix{ & Q\ar[d]^{\psi_0\omega}\ar[dl]_{\psi_1} & \\ T_1\ar[r]^{g_1}~& T_0\ar[r]& \cok g_1 .}\]}Now using the fact that $\omega$ is non-zerodivisor, we deduce that $\varphi\psi=\id_{(Q\st{\omega}\rt Q)}$, and so $(Q\st{\omega}\rt Q)$ is a projective object of $\mon(\omega, \cp )$, as required.
\end{proof}}

\begin{lemma}\label{fac}Let $(P\st{f}\rt Q)$ be an arbitrary object of $\mon(\omega, \cp)$. Then $(Q\oplus P\st{l}\rt P\oplus Q)$ with  $l=\tiny {\left[\begin{array}{ll} -f_{\si} & \id \\ {0} & {f} \end{array} \right]},$  is a projective object of $\mon(\omega, \cp)$ with $\cok l=Q/{\omega Q}$. In particula, $(P\st{f}\rt Q)$ is a homomorphic image of a projective object in $\mon(\omega, \cp)$.
\end{lemma}
\begin{proof}Since $ff_{\si}=\omega.\id_{Q}$, we may have the following commutative diagram with exact rows;\[\xymatrix{&\\ 0 \ar[r] & Q\ar[r]^{[\id~~f_{\si}]^t} \ar[d]_{\omega} & Q\oplus P\ar[r]^{[-f_{\si}~~\id]} \ar[d]_{l} & P \ar[r] \ar@{=}[d] & 0 \\  0 \ar[r] & Q \ar[r]^{[0~~\id]^t} & P\oplus Q \ar[r]^{[\id~~0]} & P \ar[r] & 0,}\]where $l=\tiny {\left[\begin{array}{ll} -f_{\si} & \id \\ {0} & {f} \end{array} \right].}$ Applying Snake lemma and using the fact that $\omega$ is non-zerodivisor, guarantee that $l$ is a monomorphism with $\cok l=Q/{\omega Q}$. Namely, $(Q\oplus P\st{l}\rt P\oplus Q)$ is an object of $\mon(\omega, \cp)$. Moreover, since by Lemma \ref{proj}, $(Q\st{\omega}\rt Q)$ and $(P\st{\id}\rt P)$ are projective objects of $\mon(\omega, \cp)$, the same is true for $(Q\oplus P\st{l}\rt P\oplus Q)$. Now the epimorphism $\pi=(\pi_1, \pi_0):(Q\oplus P\st{l}\rt P\oplus Q)\rt(P\st{f}\rt Q)$, where $\pi_1$ and $\pi_0$ are projections, completes the proof.
\end{proof}

\begin{remark}\label{con}Assume that $\psi=(\psi_1, \psi_0):(P\st{f}\rt Q)\lrt(P'\st{f'}\rt Q')$ is a morphism in $\mon(\omega, \cp)$. As mentioned in the proof of Lemma \ref{tt}, there is a short exact sequence $0\lrt(P'\st{f'}\rt Q')\st{[\id~~0]^t}\lrt(P'\oplus Q\st{c}\rt Q'\oplus P)\st{[0~~\id]}\lrt(Q\st{-f_{\si}}\rt P)\lrt 0$ in $\mon(\omega, \cp)$, where the middle term is $C(\psi)$. Then one may observe that $C([\id~~0]^t)\cong (Q\st{-f_{\si}}\rt P)$, in $\HT\mon(\omega, \cp)$. In this direction, consider the following short exact sequence in $\mon(\omega, \cp)$; $$0\lrt(P'\oplus Q'\st{l}\rt Q'\oplus P')\st{[\id~~0]^t}\lrt (P'\oplus Q\oplus Q'\st{l'}\rt Q'\oplus P\oplus P')\st{[0~~\id]}\lrt (Q\st{-f_{\si}}\rt P)\rt 0,$$where $l=\tiny {\left[\begin{array}{ll} f' & \id \\ {0} & {-f'_{\si}} \end{array} \right]}$ and the middle term of the sequence is $C([\id~~0]^t)$, and then $l'=\tiny {\left[\begin{array}{lll} f' & \psi_0 & \id\\ 0 & -f_{\si} & 0 \\ 0 & 0 & -f'_{\si}\end{array} \right].}$ It should be noted that the middle term of the latter exact sequence is the mapping cone of the injection map $(P'\st{f'}\rt Q')\st{[\id~~0]^t}\lrt (P'\oplus Q\st{c}\rt Q'\oplus P)$. In view of Lemma \ref{fac}, the left term of the above sequence is a projective object of $\mon(\omega, \cp)$. This yields that, the middle and the right terms of the sequence are isomorphic in $\HT\mon(\omega, \cp)$, as claimed. In a similar way, it can be seen that the mapping cone of the morphism $((P'\oplus Q\st{c}\rt Q'\oplus P)\st{[0~~\id]}\lrt(Q\st{-f_{\si}}\rt P))$ is isomorphic to $(P'\st{f'}\rt Q')$ in $\HT\mon(\omega, \cp)$.
\end{remark}

Now, we  prove  that $\HT\mon(\omega, \cp)$ is a triangulated category. The proof should be compared with the argument given in the proof of \cite[Theorem 6.7]{hjtriang}.

\begin{prop}\label{proptri}The triple $(\HT\mon(\omega, \cp), \si, \Delta)$ is a triangulated category.
\end{prop}
\begin{proof}We need to prove that the class of exact triangles satisfies Verdier's axioms TR1-TR4 stated in \cite[(1.1)]{happel1988triangulated}.\\ (TR1) First, one should note that, by our definition, every morphism in $\HT\mon(\omega, \cp)$ is embedded in an exact triangle and also, a triangle which is isomorphic to an exact triangle is also an exact triangle. So, we need to show that for a given object
$(P\st{f}\rt Q)\in\HT\mon(\omega, \cp)$, there exists an exact triangle $(P\st{f}\rt Q)\st{\id}\lrt(P\st{f}\rt Q)\lrt 0 \lrt(Q\st{-f_{\si}}\rt P)$. By our construction, there is an exact triangle $(P\st{f}\rt Q)\st{\id}\lrt(P\st{f}\rt Q)\st{[\id~~0]^t}\lrt C(\id)\st{[0~~\id]}\lrt(Q\st{-f_{\si}}\rt P)$, where $C(\id)=(P\oplus Q\st{c}\rt Q\oplus P)$, with $c=\tiny {\left[\begin{array}{ll} f & \id \\ {0} & {-f_{\si}} \end{array} \right].}$ By Lemma \ref{fac}, $C(\id)$ is a projective object of $\mon(\omega, \cp)$ and so, it
 is isomorphic to the zero object in the homotopy category $\HT\mon(\omega, \cp)$, as needed.\\
(TR2) We have to show that any rotation of an exact triangle in $\HT\mon(\omega, \cp)$ is also exact. Without loss of generality, one may consider standard triangles. So, for a given standard triangle $(P\st{f}\rt Q)\st{\psi}\lrt (P'\st{f'}\rt Q')\st{[\id~~0]^t}\lrt C(\psi)\st{[0~~\id]}\lrt(Q\st{-f_{\si}}\rt P)$, we shall prove that the triangle $(P'\st{f'}\rt Q')\st{[\id~~0]^t}\lrt C(\psi)\st{[0~~\id]}\lrt(Q\st{-f_{\si}}\rt P)\st{-\si\psi}\lrt(Q'\st{-f'_{\si}}\rt P')$ is an exact triangle. According to Remark \ref{con}, $(Q\st{-f_{\si}}\rt P)$ isomorphic to $C([\id~~0]^t)$. Consequently, the latter triangle is
 isomorphic in $\HT\mon(\omega, \cp)$ to the standard triangle $(P'\st{f'}\rt Q')\st{[\id~~0]^t}\lrt C(\psi)\st{}\lrt C([\id~~0]^t)\st{}\lrt(Q'\st{-f'_{\si}}\rt P')$, giving the desired result.\\
(TR3) Again it would be enough to verify the axiom for standard triangles. So, assume  that
{\footnotesize \[\xymatrix{(P\st{f}\rt Q)\ ~\ar[r]^{\psi}  \  \ \ar[d]_{\gamma}& \ \  (P'\st{f'}\rt Q')\ar[r]^{[\id~~0]^t} \ \ \ar[d]_{\epsilon}&  \ \ C(\psi)\ar[r]^{[0~~\id]} \ \  & \  (Q\st{-f_{\si}}\rt P)\ar[d]_{\si\gamma}\\ (P_1\st{g}\rt Q_1)\ ~\ar[r]^{\epsilon'} \ \ &\ \ (P'_1\st{g'}\rt Q'_1)\ar[r]^{[\id~~0]^t} \ \ & \ \ C(\epsilon')\ar[r]^{[0~~\id]} \ \ & (Q_1\st{-g_{\si}}\rt P_1),}\]}is a diagram in $\HT\mon(\omega, \cp)$ such that the left square is commutative. We must find a morphism $\eta:C(\psi)\lrt C(\epsilon')$ which makes the diagram commute. By our hypothesis, there are homotopy morphisms $s_0:Q\rt P'_1$ and $s_1:P\rt Q'_1$ such that $(\epsilon_0\psi_0-\epsilon'_0\gamma_0)f-g's_0f=\omega.s_1$. Thus, in order to complete the diagram to a morphism of triangles,  one may define $\eta=(\eta_1, \eta_0):C(\psi)\lrt C(\epsilon')$ by letting $\eta_1=\tiny {\left[\begin{array}{ll} \epsilon_1 & s_0 \\ 0 & \gamma_0 \end{array} \right]}$ and $\eta_0=\tiny {\left[\begin{array}{ll} \epsilon_0 & s_1 \\ 0 & \gamma_1 \end{array} \right].}$ Now it is easily checked that the constructed squares commute.\\(TR4) As the previous, it suffices to check the octahedral axiom for standard triangles. Assume that $(P\st{f}\rt Q)\st{\psi}\lrt (P'\st{f'}\rt Q')\st{\eta}\lrt (P''\st{f''}\rt Q'')$ is a composition of morphisms in $\HT\mon(\omega, \cp)$. We must prove that there exists the following commutative diagram;
 \[\xymatrix{(P\st{f}\rt Q)~\ar[r]^{\psi}\ar@{=}[d]& (P'\st{f'}\rt Q')\ar[r]^{[\id~~0]^t}\ar[d]_{\eta}& C(\psi)\ar[r]^{[0~~\id]}\ar[d]_{\varphi} &(Q\st{-f_{\si}}\rt P)\ar@{=}[d]\\ (P\st{f}\rt Q)~\ar[r]^{\eta\psi}\ar[d]_{\psi}& (P''\st{f''}\rt Q'')\ar[r]^{[\id~~0]^t}\ar@{=}[d]& C(\eta\psi)\ar[r]^{[0~~\id]}\ar[d]_{\gamma} &(Q\st{-f_{\si}}\rt P)\ar[d]_{\si(\psi)}\\ (P'\st{f'}\rt Q')~\ar[r]^{\eta}\ar[d]_{[\id~~0]^t}& (P''\st{f''}\rt Q'')\ar[r]^{[\id~~0]^t}\ar[d]_{[\id~~0]^t}& C(\eta)\ar[r]^{[0~~\id]}\ar@{=}[d] &(Q'\st{-f'_{\si}}\rt P')\ar[d]_{\si([\id~~0]^t)}\\  C(\psi)\ar[r]^{\varphi} & C(\eta\psi)\ar[r]^{\gamma} & C(\eta)\ar[r]^{\delta}&\si C(\psi),}\]in $\HT\mon(\omega, \cp)$ such that  rows are exact triangles. To this end, by applying TR1 and TR3, we only need to show the bottom row is an exact triangle. This will be done by proving that it is isomorphic to the standard triangle $C(\psi)\st{\varphi}\lrt C(\eta\psi)\st{[\id~~0]^t}\lrt C(\varphi)\st{[0~~\id]}\lrt\si C(\psi)$. In order to construct an isomorphism between these triangles, we may take the identity morphisms for the first, second and fourth entries, and for the third entry, we define the following morphism; $$\epsilon=(\epsilon_1, \epsilon_0):C(\eta)=(P''\oplus Q'\st{l}\rt Q''\oplus P')\lrt C(\varphi)=(P''\oplus Q\oplus Q'\oplus P\st{l'}\rt Q''\oplus P\oplus P'\oplus Q),$$ by setting $\epsilon_i=\tiny {\left[\begin{array}{llll} \id & 0 &0 &0 \\ 0 & 0& \id& 0 \end{array} \right]^t}$, where $l=\tiny {\left[\begin{array}{ll} f'' & \eta_0 \\ 0 & -f'_{\si} \end{array} \right]}$ and $l'=\tiny {\left[\begin{array}{llll} f'' & \eta_0\psi_0 &\eta_0 &0 \\ 0 & -f_{\si}& 0& \id \\ 0 & 0& f'_{\si}& \psi_1\\ 0 & 0& 0& -f\end{array} \right]}$.  Since $\gamma=(\gamma_1,\gamma_0)$ and $\delta=(\delta_1, \delta_0)$ are morphisms with $\gamma_1=\tiny {\left[\begin{array}{ll} \id & 0 \\ 0 & \psi_0 \end{array} \right]}$, $\gamma_0=\tiny {\left[\begin{array}{ll} \id & 0 \\ 0 & \psi_1 \end{array} \right]}$ and $\delta_i=\tiny {\left[\begin{array}{ll} 0 & \id \\ 0 & 0 \end{array} \right]}$, we infer that $[0~~\id]\epsilon=\delta$ and $\epsilon\gamma-[\id~~0]^t$ is null-homotopic with the homotopy morphisms $s_0=\tiny {\left[\begin{array}{llll} 0 & 0 &0 &0 \\ 0 & 0& 0& -\id \end{array} \right]^t}$ and $s_1=\tiny {\left[\begin{array}{llll} 0 & 0 &0 &0 \\ 0 & 0& 0& \id \end{array} \right]^t}$. So, it remains to show that $\epsilon$ is an isomorphism in $\HT\mon(\omega, \cp)$. To this end, consider the short exact sequence, $0\lrt C(\eta)\st{\epsilon}\lrt C(\varphi)\lrt (Q\oplus P\st{l''}\rt P\oplus Q)\lrt 0$ in $\mon(\omega, \cp)$, with $l''=\tiny {\left[\begin{array}{ll} f_{\si} & \id \\ 0 & -f \end{array} \right]}$. According to Lemma \ref{fac}, the right term is a projective object of $\mon(\omega, \cp)$, and so, $\epsilon$ will be an isomorphism in $\HT\mon(\omega, \cp)$. So the proof is completed.
\end{proof}

\begin{s}{\sc Gorenstein projective modules.} An acyclic complex of projective $\Lambda$-modules; $\mathbf{P}_{\bullet}:\cdots\lrt P_{n+1}\st{d_{n+1}}\lrt P_n\st{d_n}\lrt P_{n-1}\st{d_{n-1}}\lrt\cdots$ is called {\em totally acyclic}, if the acyclicity is preserved by $\Hom_{\Lambda}(-, P )$ for every projective $\Lambda$-module $P$.  A $\Lambda$-module $M$ is said to be {\em Gorenstein projective}, if it is a syzygy of a totally acyclic complex of projective modules. Clearly, every projective module is Gorenstein projective.

 It is known that over a Gorenstein ring $\Lambda$, every acyclic complex is totally acyclic, and also $d$-th syzygy of any $\Lambda$-module is Gorenstein projective, where $d=\dim \Lambda$; see \cite[Theorem 10.2.14]{enochs2011relative}. Finitely generated Gorenstein projective modules over a noetherian ring are introduced by Auslander and Bridger under the name ``modules of G-dimension zero'' \cite{auslander1969stable}.  Over a commutative Gorenstein ring, these modules are equal to the maximal Cohen-Macaulay modules. The category of all (finitely generated) Gorenstein projective $\Lambda$-modules, will be depicted by $\gp(\Lambda)$.
\end{s}

\begin{remark}\label{rr}Assume that $(P\st{f}\rt Q)$ is an arbitrary object of $\mon(\omega, \cp)$. Since $ff_{\si}=\omega.\id_Q$ and $f_{\si}f=\omega.\id_P$, one may get an exact sequence of $R$-modules; $0\rt\cok f\rt P/{\omega P}\st{\bar{f}}\rt Q/{\omega Q}\rt\cok f\rt 0$, and so, $\cok f$ has a 2-periodic projective resolution. We should emphasize that, as $\cok f$ is an $R$-module, the equality $\cok f=\cok\bar{f}$ holds. Let us explain the equality $\ker\bar{f}=\cok f$. Applying the functor $-\otimes_{S}R$ to the short exact sequence of $S$-modules $0\rt P\st{f}\rt Q\rt\cok f\rt 0$, gives rise to the exact sequence $0\rt\Tor_1^S(\cok f, R)\rt P/{\omega P}\st{\bar{f}}\rt Q/{\omega Q}\rt \cok f\rt 0$. So it suffices to show that $\Tor_1^S(\cok f, R)=\cok f$. To derive this, one may apply the functor $\cok f\otimes_{S}-$ to the short exact sequence of $S$-modules; $0\rt S\st{\omega}\rt S\rt R\rt 0$ and use the fact that $\cok f$ is annihilated by $\omega$. Similarly, we have also an exact sequence of $R$-modules; $0\rt\cok f_{\si}\rt Q/{\omega Q}\st{\bar{f}_{\si}}\rt P/{\omega P}\rt\cok f_{\si}\rt 0$. It should be noted that these exact sequences yield that $\syz_{{R}}(\cok f)=\cok f_{\si}$ and $\syz_{{R}}(\cok f_{\si})=\cok f$. In particular, we get the acyclic complex of projective $R$-modules; $\cdots\rt P/{\omega P}\st{\bar{f}}\rt Q/{\omega Q}\st{\bar{f}_{\si}}\rt P/{\omega P}\st{\bar{f}}\rt Q/{\omega Q}\st{\bar{f}_{\si}}\rt\cdots$. Since the projective dimensions of $\cok f$ and $\cok f_{\si}$ over $S$ are at most one, applying \cite[Lemma 2(i), page 140]{mat} enables us to deduce that $\Ext^i_R(\cok f, R)=0=\Ext^i_R(\cok f_{\si}, R)$ for all $i\geq 1$. Consequently, the latter acyclic complex is totally acyclic, and so,  $\cok f$ is a  Gorenstein projective $R$-module; see also \cite[Lemma 3.1]{bergh2015complete}.
\end{remark}

The result below, which is the main result of this section, clarifies the connection between the homotopy category $\HT\mon(\omega, \cp)$ and the stable category of Gorenstein projective $R$-modules $\ugp(R)$.

\begin{theorem}\label{ds}There is a fully faithful functor $T:\HT\mon(\omega, \cp)\lrt\ugp(R)$, sending each object $(P\st{f}\rt Q)$ to $\cok f$.
\end{theorem}
\begin{proof}Assume that $(P\st{f}\rt Q)$ is an arbitrary object of $\HT\mon(\omega, \cp)$. As noted in Remark \ref{rr}, $\cok f$ is a Gorenstein projective $R$-module. Moreover, any morphism $\psi=(\psi_1, \psi_0):(P\st{f}\rt Q)\lrt(P'\st{f'}\rt Q')$ gives us a unique $S$ (and also $R$)-homomorphism $h:\cok f\rt\cok f'$. Now, we set $T(\psi):=h$. Assume that $\psi=(\psi_1, \psi_0)$ is null-homotopic. We must show that $T(\psi)$ is zero in $\ugp(R)$. Take $S$-homomorphisms $s_0:Q\rt P'$ and $s_1:P\rt Q'$ such that $f'\psi_1-f's_0f=\omega. s_1$. This in conjunction with Lemma \ref{fac}, gives us  the following commutative diagram with exact rows;
$$\xymatrix{0\ar[r] & P\ar[r]^{f} \ar[d]_{[s_1~~\psi_1]^t} &\ Q \ar[d]_{[s_0~~\psi_0]^t}\ar[r] & \cok f \ar[d]_{\alpha}\ar[r]& 0\\ 0\ar[r]& Q'\oplus P' \ar[r]^{l}\ar[d]_{[0~~\id]} & P'\oplus Q' \ar[r] \ar[d]^{[0~~\id]}  & Q'/{\omega Q'} \ar[d]_{\beta}\ar[r]&0 \\0\ar[r]&P' \ar[r]^{f'} &  Q' \ar[r] & \cok f'\ar[r]& 0,}$$where $l=\tiny {\left[\begin{array}{ll} -f'_{\si} & \id \\ {0} & {f'} \end{array} \right].}$ Since the compositions of the left and the middle columns are $\psi_1$ and $\psi_0$, respectively,
$T(\psi)=\beta\alpha$ factors through the projective $R$-module $Q'/{\omega Q'}$. Thus it is zero in  $\ugp(R)$.  Consequently, $T$ is well-defined.\\
The functor $T$ is full: assume that $(P\st{f}\rt Q)$ and $(P'\st{f'}\rt Q')$ are two objects of $\HT\mon(\omega, \cp)$, and a morphism $h:\cok f\rt\cok f'$ is given. Since $Q$ is a projective $S$-module, one may get the following commutative diagram with exact rows; $$\xymatrix{0\ar[r] & P\ar[r]^{f} \ar[d]_{\psi_1} &\ Q \ar[d]_{\psi_0}\ar[r] & \cok f \ar[d]_{h}\ar[r]& 0 \\0\ar[r]&P' \ar[r]^{f'} &  Q' \ar[r] & \cok f'\ar[r]& 0.}$$ This means that $\psi=(\psi_1,\psi_0):(P\st{f}\rt Q)\lrt(P'\st{f'}\rt Q')$ is a morphism in $\HT\mon(\omega, \cp)$ and $T(\psi)=h$.\\
Now we prove that the functor $T$ is faithful. To do this, assume that $\psi=(\psi_1, \psi_0):(P\st{f}\rt Q)\lrt(P'\st{f'}\rt Q')$ is a morphism in $\HT\mon(\omega, \cp)$ such that $T(\psi)=h=0$ in $\ugp(R)$. We have to show that $\psi=0$. Since $h=0$ in $\ugp(R)$, it factors through a projective $R$-module $P_1/{\omega P_1}$, for some projective $S$-module $P_1$, because $R$ is the local  factor ring $S/{(\omega)}$. Hence, one may obtain the following commutative diagram with exact rows; $$\xymatrix{0\ar[r] & P\ar[r]^{f} \ar[d]_{s} &\ Q \ar[d]_{t}\ar[r] & \cok f \ar[d]_{\alpha}\ar[r]& 0\\ 0\ar[r]& P_1\ar[r]^{\omega}\ar[d]_{s'} & P_1 \ar[r] \ar[d]_{t'}  & P_1/{\omega P_1} \ar[d]_{\beta}\ar[r]&0 \\0\ar[r]&P' \ar[r]^{f'} &  Q' \ar[r] & \cok f'\ar[r]& 0.}$$It should be noted that the existence of morphisms $t$ and $t'$ come from the projectivity of $Q$ and $P_1$. So, it is routine to check that there is an $S$-homomorphism $s_0:Q\rt P'$ such that $f's_0=\psi_0-t't$. Thus, one gets the equality $\psi_0f-f's_0f=t'tf$. Now setting $s_1:=t's$, we have the equality $\omega. s_1=t'tf$. Combining this with the former equality, gives us the equality $\psi_0f-f's_0f=\omega.s_1$, meaning that $\psi=(\psi_1, \psi_0)$ is null-homotopic, and then, it is zero in $\HT\mon(\omega, \cp)$. Thus the proof is completed.
\end{proof}

It is known that the category $\gp(R)$  is a Frobenius category. So the stable category $\ugp(R)$ admits a natural structure of a triangulated category with the quasi-inverse of the syzygy functor $\syz^{-1}:\ugp(R)\rt\ugp(R)$ as suspension. Indeed, cosyzygies with respect to injective objects in the Frobenius category $\gp(R)$ is taken; see \cite[Chapter I, Section 2]{happel1988triangulated} for more details.

\begin{prop}With the notation above, $T$ is a triangle functor.
\end{prop}
\begin{proof}
\new{}
It suffices to show that for any object ${\bf f}=(P\st{f}\rt Q)\in\HT\mon(\omega, \cp)$, $T\si({\bf f})=\syz^{-1}T({\bf f})$. By our definition,  $T\si({\bf f})=T({-{\bf f}_{\si}})=\cok(-f_{\si})$, where $-{\bf f}_{\si}=(Q\st{-f_{\si}}\rt P)$. Moreover, $\syz^{-1}T({\bf f})=\syz^{-1}(\cok f)$. So we need to examine that $\syz^{-1}(\cok f)=\cok(-f_{\si})$. To do this, consider the following commutative diagram; $${\fo\xymatrix{&& 0\ar[d]& 0\ar[d]\\ 0\ar[r] & P\ar[r]^{f} \ar[d]_{-\id} &\ Q \ar[d]_{-f_{\si}}\ar[r] & \cok f \ar[d]\ar[r]& 0 \\0\ar[r]&P \ar[r]^{\omega} & P \ar[r]\ar[d] & P/{\omega P}\ar[r]\ar[d]& 0.\\ && \cok(-f_{\si})\ar[r]^{\id}\ar[d] & \cok(-f_{\si})\ar[d]\\ && 0&0}}$$Now, according to the right column, we have that $\cok(-f_{\si})=\syz^{-1}(\cok f)$, giving the desired result.
\end{proof}

The singularity category $\ds(R)$  is by definition the Verdier quotient of the bounded derived category $\db(R)$ of $R$ by the perfect complexes. This category measures the homological singularity of $R$ in the sense that $R$ has finite global dimension if and only if its singularity category is trivial. This notion was introduced by Buchweitz \cite{buchweitz1987maximal} in the 1980s, and studied actively ever since the relation with mirror symmetry was found by Orlov \cite{orlov2003triangulated}.

It is known that the functor $T':\ugp(R)\lrt\ds(R)$, sending each object to its stalk complex, is fully faithful; see \cite[Theorem 3.1] {bergh2015gorenstein}. By gluing this with Theorem \ref{ds}, we have the result below.

\begin{cor}\label{ccc}There is a fully faithful functor $F:\HT\mon(\omega, \cp)\lrt\ds(R)$, sending each object $(P\st{f}\rt Q)$ to $\cok f$, viewed as a stalk complex.
\end{cor}

\section{Specifying the Auslander-Reiten translation in $\mon(\omega, \cp)$}
This section aims to determine the Auslander-Reiten translation in the category $\mon(\omega, \cp)$. To be precise, assume that the factor ring $R$ is a complete Gorenstein ring with $d=\di R$.  It is known that {for any non-projective indecomposable Gorenstein projective $R$-module $M$, there is an almost split sequence ending at $M$} in the category $\gp(R)$, provided that $R$ is an isolated singularity. In particular, the Auslander-Reiten translation is given by $\tau(-)=\Hom_R(\syz^d\Tr_R(-), R)$, where $\Tr_R(-)$ stands for the Auslander transpose and  $\syz^i(-)$ is the $i$-th syzygy functor defined as usual by taking the kernels of the projective covers consecutively. By using Theorem \ref{ds}, we give an explicit description of the Auslander-Reiten translation in the category $\mon(\omega, \cp)$. Moreover, it is proved that {each indecomposable non-projective object of $\mon(\omega, \cp)$ appears as the right term of an almost split sequence in $\mon(\omega, \cp)$}, whenever $R$ is an isolated singularity. We begin by recalling the definition of almost split sequences.

\begin{s}{\sc Almost split sequences.} Let $\mathcal{A}$ be an exact category. A morphism
$f:B\rt C$  in $\mathcal{A}$ is said to be right almost split provided it is not a split epimorphism and
every morphism $h:X\rt C$ which is not a split epimorphism factors through $f$. {Dually, a morphism
$g: A\rt B$  is left almost split, if it is not a split monomorphism and each morphism $h:A\rt C$ which
is not a split monomorphism factors through $g$. A conflation  $\eta:A\st{f}\rt B\st{g}\rt C$ is called
almost split, provided that $f$ is left almost split and ${g}$ is right almost split; see \cite{ji}.
We remark that, since this sequence is unique up to isomorphism for $A$ and for $C$,
we may write $A=\tau_{\A}{C}$ and $C=\tau_{\A}^{-1}A$.
It is known that if $\eta:A\st{f}\rt B\st{g}\rt C$ is an almost split conflation, then $A$ and $C$ have local endomorphism ring; see \cite[Proposition II.4.4]{auslander1967functors}.}
 In the reminder, we fix the notation $\tau:=\tau_{\gp(R)}$.

The existence of almost split sequences also known as the Auslander-Reiten sequences is a fundamental and important component in the
study of Auslander-Reiten theory, which was introduced by Auslander and Reiten in \cite{auslander1975representation}, where they have proved
the first existence theorem for the category of finitely generated modules over an artin algebra. This theory rapidly developed in various contexts such as orders over Gorenstein rings \cite{auslander1967functors}
and the category of maximal Cohen-Macaulay modules over a Henselian Cohen-Macaulay local ring which admits a canonical module; see \cite{leuschke2012cohen,yoshino1990cohen}.
For the terminology and background on almost split morphisms, we refer the reader to  \cite{auslander1995representation, ausm, liukrull}. Moreover,  the reader may consult \cite{inp} for a general setting of the Auslander-Reiten theory.
\end{s}

The next result gives the precise determination of the Auslander-Reiten translation of some specific objects in the category of $\gp(R)$.

\begin{prop}\label{dim}Let $M$ be an indecomposable non-projective $R$-module which has a 2-periodic projective resulotion. Then $\tau(M)=\syz M$, if $d=\di R$ is an odd number, and $\tau(M)=M$, whenever $\di R$ is even.
\end{prop}
\begin{proof}By our assumption, there exists an exact sequence of $R$-modules; $0\rt M\rt P_1\rt P_0\rt M\rt 0$ in which $P_0,P_1\in\cp(R)$. Since $R$ is Gorenstein, { every acyclic complex of projectives is totally acyclic and so,}  the sequence $0\rt M^*\rt P_0^*\rt P_1^*\rt M^*\rt 0$ is also exact, where $(-)^*=\Hom_R(-, R)$; see \cite[Corollary 5.5]{ik}. So $\Tr(M)=M^*$, and in particular, $\Tr(M)$ admits a 2-periodic projective resolution, because $P_0^*, P_1^*\in\cp(R)$. Consequently, for any even integer $n$, $\syz^n\Tr(M)=\Tr(M)$, and then $(\syz^n\Tr(M))^*=M^{**}=M$. Next assume that $n$ is odd. So $\syz^n\Tr(M)=\syz(\syz^{n-1}\Tr(M))=\syz\Tr(M)$. Now considering the short exact sequence of $R$-modules; $0\rt(\syz\Tr(M))^*\rt P_0\rt M\rt 0$, we infer that $(\syz\Tr(M))^*=\syz M$. Finally, since $\tau(M)=(\syz^d\Tr(M))^*$, one obtains that $\tau(M)=M$, whenever $d$ is even, and $\tau(M)=\syz M$, if $d$ is odd. Thus the proof is finished.
\end{proof}

\begin{prop}\label{rrrr}A morphism $\psi=(\psi_1, \psi_0):(P\st{f}\rt Q)\lrt(P'\st{f'}\rt Q')$ is  null-homotopic if and only if  $\psi$ factors through a projective object of $\mon(\omega, \cp)$.
\end{prop}
\begin{proof}
First, we deal with the `only if' part. Since $\psi=(\psi_1, \psi_0):(P\st{f}\rt Q)\lrt(P'\st{f'}\rt Q')$ is a null-homotopic morphism, there are $S$-homomorphisms $s_0:Q\rt P'$ and $s_1:P\rt Q'$ such that $f'\psi_1-f's_0f=\omega. s_1$. So, one may obtain the following commutative diagram; $$\xymatrix{ P\ar[r]^{f} \ar[d]_{[s_1~~\psi_1]^t} &\ Q \ar[d]_{[s_0~~\psi_0]^t}\\ Q'\oplus P' \ar[r]^{l}\ar[d]_{[0~~\id]} & P'\oplus Q' \ar[d]^{[0~~\id]}  \\ P' \ar[r]^{f'} &  Q',}$$where $l=\tiny {\left[\begin{array}{ll} -f'_{\si} & \id \\ {0} & {f'} \end{array} \right].}$ That is, the morphism $\psi$ factors through $(Q'\oplus P')\st{l}\rt (P'\oplus Q')$ which is a projective object of $\mon(\omega, \cp)$, because of Lemma \ref{fac}. For the `if' part, at first we should highlight that, according to Lemma \ref{fac},  each projective object in $\mon(\omega, \cp)$ is equal to direct summands of finite direct sums of an object of the form $(P_1\oplus Q_1\st{l}\rt P_1\oplus Q_1)$, where $P_1, Q_1\in\cp(S)$ and $l=\tiny {\left[\begin{array}{ll} \omega & 0 \\ {0} & {\id} \end{array} \right].}$ This fact  allows us to  assume that the morphism $\psi=(\psi_1, \psi_0):(P\st{f}\rt Q)\lrt(P'\st{f'}\rt Q')$ factors through a projective object $(P_1\oplus Q_1\st{\omega\oplus\id}\lrt P_1\oplus Q_1)$. Namely, we have the following commutative diagram;
$$\xymatrix{ P\ar[r]^{f} \ar[d]_{[\alpha_1~~\beta_1]^t} &\ Q \ar[d]_{[\alpha_0~~\beta_0]^t}\\ P_1\oplus Q_1 \ar[r]^{\omega\oplus\id}\ar[d]_{[\alpha'_1~~\beta'_1]} & P_1\oplus Q_1 \ar[d]^{[\alpha'_0~~\beta'_0]}  \\ P' \ar[r]^{f'} &  Q',}$$
such that the composition of morphisms in the left (resp. right) column is $\psi_1$ (resp. $\psi_0$).  Now set $s_0:={[\alpha'_1~~\beta'_1]}{[\alpha_0~~\beta_0]^t}-f'_{\si}\psi_0$ and $s_1:={[\alpha'_0~~\beta'_0]}{[\alpha_1~~\beta_1]^t}$. Hence, using the facts that $\omega.\alpha_1=\alpha_0f$, $\beta_0f=\beta_1$, $\omega.\alpha'_0=f'\alpha'_1$ and $f'\beta'_1=\beta'_0$, one deduce that  $f'\psi_1-f's_0f=\omega.s_1$, meaning that $\psi=(\psi_1, \psi_0)$ is a null-homotopic morphism. Thus the proof is completed.
\end{proof}

\begin{prop}\label{p1}The following assertions are satisfied: \\(1) A given object $(P\st{f}\rt Q)\in\mon(\omega, \cp)$ is indecomposable if and only if $(Q\st{f_{\si}}\rt P)$ is so. \\(2) Let {$(P\st{f}\rt Q)$ be an indecomposable object of $\mon(\omega, \cp)$. Then $\cok f$ is an indecomposable object of $\gp(R)$.} \\ (3) Let $(P\st{f}\rt Q)$  be an object of $\mon(\omega, \cp)$ such that $Q\st{\pi}\rt\cok f$ is the projective cover of  $\cok f$ over $S$. If $\cok f$ is an indecomposable non-projective $R$-module, then $(P\st{f}\rt Q)$ is an indecomposable non-projective object of  $\mon(\omega, \cp)$.\\ (4) {If $(P\st{f}\rt Q)$ is an indecomposable non-projective object of $\mon(\omega, \cp)$, then $Q\rt\cok f$ is the projective cover of $\cok f$ over $S$. \\ (5) If $(P\st{f}\rt Q)$ is an indecomposable non-projective object of $\mon(\omega, \cp)$, then it has local endomorphism ring.}
\end{prop}
\begin{proof}(1) This comes up directly from the equality $(f_{\si})_{\si}=f$.\\
(2) Assume on the contrary that $\cok f$ is descomposable, and so, we may write $\cok f=X_1\oplus X_2$. As $S$ is a local ring, applying \cite[Theorem 5.3.3]{enochs2011relative} gives us the projective covers $Q_1\st{\pi_1}\rt X_1$ and $Q_2\st{\pi_2}\rt X_2$. Thus $Q_1\oplus Q_2\st{\pi_1\oplus \pi_2}\rt X_1\oplus X_2$ will be the projective cover of $X_1\oplus X_2$ over $S$; see \cite[Corollary 5.5.2]{enochs2011relative}. Since $X_1$ and $X_2$ have projective dimension at most one over $S$, we may take a short exact sequence of $S$-modules; $0\rt P_1\oplus P_2\rt Q_1\oplus Q_2\rt \cok f\rt 0$ in which $P_1, P_2\in\cp(S)$. In particular, we get the following commutative diagram with exact rows;
$$\xymatrix{0\ar[r] & P_1\oplus P_2\ar[r] \ar[d]_{g} &\ Q_1\oplus Q_2 \ar[d]_{h}\ar[r] & \cok f \ar@{=}[d]\ar[r]& 0 \\0\ar[r]&P \ar[r]^{f} &  Q \ar[r] & \cok f\ar[r]& 0,}$${where $h$ is a split monomorphism, because $Q_1\oplus Q_2\rt\cok f$ is the projective cover, and so, $Q_1\oplus Q_2$ will be a direct summand of Q. Thus $\cok h=T$ is a projective $S$-module, implying that $g$ is also a split monomorphism, thanks to the fact that $\cok g=T$. Hence, $\cok(g, h)=(T\st{id}\rt T)$ is a projective object of  $\mon(\omega, \cp)$.} This, particularly means that $(P_1\rt Q_1)$ is  a non-zero direct summand of the indecomposable object $(P\st{f}\rt Q)$, which is impossible. Consequently, $\cok f$ is an indecomposable object of $\gp(R)$.\\
(3) Assume to the contrary that $(P\st{f}\rt Q)$ is a decomposable object of $\mon(\omega, \cp)$. {Since $\cok f$ is indecomposable and $R$ is complete, $\End_R(\cok f)$ will be a local ring, and so, the same will be true for $\uEnd_R(\cok f)$. Moreover, in view of Theorem \ref{ds}, we have an isomorphism $\uEnd_R(\cok f)\cong\uEnd_{\mon}(P\st{f}\rt Q)$, implying that the latter is also a local ring. This, in conjunction with $(P\st{f}\rt Q)$ being decomposable, forces $(P\st{f}\rt Q)$ to have an indecomposable projective direct summand $(Q'\rt Q')$. According to Lemma \ref{proj}, this should be either $(Q'\st{\id}\rt Q')$ or $(Q'\st{\omega}\rt Q')$. If the latter one takes place, then $Q'/{\omega Q'}$ will be a (projective) direct summand of $\cok f$, which is impossible. Consequently, we have $(P\st{f}\rt Q)=(Q_1\st{g}\rt P_1)\oplus (Q'\st{\id}\rt Q')$, and then, {by applying Lemma \ref{lem1} we have} $(Q\st{f_{\si}}\rt P)=(P_1\st{g_{\si}}\rt Q_1)\oplus (Q'\st{\omega}\rt Q')$. This contradicts with the fact that $\cok f_{\si}=\syz_R(\cok f)$ is indecomposable, because $Q\st{\pi}\rt\cok f$ is the projective cover of  $\cok f$. Thus $(P\st{f}\rt Q)$ will be an indecomposable object of $\mon(\omega, \cp)$, and so,  the proof is completed}.\\{(4) First one should note that since $(P\st{f}\rt Q)$ is non-projective, $\cok f$ will be a non-zero $R$-module.  As $S$ is a local ring,  one may take an exact sequence of $S$-modules; $0\rt P'\rt Q'\rt\cok f\rt 0$, where $P', Q'\in\cp(S)$ and $Q'\rt\cok f$ is the projective cover of $\cok f$. So, as we have observed in the proof of the second assertion, $(P'\rt Q')$ is a direct summand of $(P\rt Q)$, and then, they will be equal, because $(P\st{f}\rt Q)$ is indecomposable.  Thus $Q\rt\cok f$ is the projective cover of $\cok f$, as needed.\\ (5) First it should be noted that for a given object $\psi=(\psi_1, \psi_0)\in\End_{\mon}(P\st{f}\rt Q)$, one may get the following commutative diagram of $S$-modules with exact rows; $$\xymatrix{0\ar[r] & P\ar[r]^{f} \ar[d]_{\psi_1} &\ Q \ar[d]_{\psi_0}\ar[r] & \cok f \ar[d]_{h}\ar[r]& 0 \\0\ar[r]&P \ar[r]^{f} &  Q \ar[r] & \cok f\ar[r]& 0.}$$Since $(P\st{f}\rt Q)$ is indecomposable, the fourth assertion yields that $Q\rt\cok f$ is the projective cover of $\cok f$, implying that $\psi=(\psi_1, \psi_0)$ is an isomorphism if and only if   $h$ is an isomorphism as $S$ (and also $R$)-homomorphism. Now assume that $\psi, \psi'\in\End_{\mon}(P\st{f}\rt Q)$ which are not isomorphism, and $h, h'$ are corresponding morphisms in $\End_R(\cok f)$. As $\psi$ and $\psi'$ are non-isomorphisms, as already observed, $h$ and $h'$ are so. According to the second assertion, $\cok f$ is indecomposable, and so, it has local endomorphism ring, thanks to the completeness of $R$.
 Consequently, $h+h'$ is a non-isomorphism, and then, the same will be true for $\psi+\psi'$, meaning that $\End_{\mon}(P\st{f}\rt Q)$ is a local ring. Hence the proof is finished.
}
\end{proof}

The next elementary result is needed for the subsequent theorem.
\begin{lemma}\label{split}Let $\A$ be an additive category with projective objects and let $g:M\rt M''$ be an epimorphism in $\A$. If there is a morphism $f:M''\rt M$ such that $\id_{M''}-gf$ factors through a projective object in $\A$, then $g$ is a split epimorphism.
\end{lemma}
\begin{proof}By our hypothesis, there are morphisms $M''\st{h_1}\rt P\st{h_2}\rt M''$ with $P$ projective in $\A$, such that $\id_{M''}-gf=h_2h_1$. Take a morphism $u:P\rt M$ in which $gu=h_2$. Set $\psi:=uh_1+f$.  So one may have the equalities; $g\psi=g(uh_1+f)=guh_1+gf=h_2h_1+gf=\id_{M''}$. This indeed means that $g$ is a split epimorphism. So the proof is finished.
 \end{proof}

\begin{theorem}\label{ass}Let $(P\st{f}\rt Q)$ be an indecomposable non-projective object of $\mon(\omega, \cp)$ and set $M=\cok f$. Then there is an almost split sequence ending at $(P\st{f}\rt Q)$ in $\mon(\omega, \cp)$ if and only if there is an almost split sequence ending at $M$ in $\gp(R)$. In particular, $\tau_{\mon}((P\st{f}\rt Q))=(P\st{f}\rt Q)$, if $\di R$ is even, and $\tau_{\mon}((P\st{f}\rt Q))=(Q\st{f_{\si}}\rt P)$, if $\di R$ is odd.
\end{theorem}
\begin{proof}First one should note that, as declared in  Remark \ref{rr}, $M$ is a Gorenstein projective $R$-module  which has a 2-periodic projective resolution. Assume that there is an almost split sequence; $0\rt\tau(M)\rt X\rt M\rt 0$ in $\gp(R)$. We would like to show that there is an almost split sequence in $\mon(\omega, \cp)$ ending at $(P\st{f}\rt Q)$. In the route to this goal, we consider two cases;\\ Case 1: Suppose that $\di R$ is odd. {Since by  Remark \ref{rr}, $M$ is a Gorenstein projective $R$-module  having a 2-periodic projective resolution, applying  Proposition \ref{dim} yields that $\tau(M)=\syz M$.} As $M$ corresponds to $(P\st{f}\rt Q)$, considering Remark \ref{rr}, $\syz M$ will correspond to the object $(Q\st{f_{\si}}\rt P)$ in $\HT\mon(\omega, \cp)$. Thus, one may get the following commutative diagram with exact rows and columns; $$\xymatrix{&0\ar[d]&0\ar[d]& 0\ar[d]&\\ 0\ar[r] & Q\ar[r] \ar[d]_{f_{\si}} &\ Q\oplus P \ar[d]\ar[r]^{[0~~\id]} & P \ar[d]_{f}\ar[r]& 0\\ 0\ar[r]& P\ar[r]^{[\id~~0]^t}\ar[d]_{\varphi_1} & P\oplus Q \ar[r]^{[0~~\id]} \ar[d]_{[f_1\varphi_1~~\beta]}  & Q \ar[d]_{\varphi_0}\ar[r]&0 \\0\ar[r]&\syz M\ar[d] \ar[r]^{f_1} &  X \ar[r]^{f_2}\ar[d] & M\ar[r]\ar[d]& 0,\\ &0&0&0&}$$where $\beta:Q\rt X$ is a morphism such that $f_2\beta=\varphi_0$. It should be remarked that $(Q\oplus P\rt P\oplus Q)$ lies in $\mon(\omega, \cp)$, because $X$ is an $R$-module. In particular, we obtain the short exact sequence, $\eta:0\lrt(Q\st{f_{\si}}\rt P)\st{\theta}\lrt(Q\oplus P\st{q}\rt P\oplus Q)\st{g}\lrt(P\st{f}\rt Q)\lrt 0$ in $\mon(\omega, \cp)$. Here $g=(g_1, g_0)$ with $g_0=[0~~\id_Q]$ and $g_1=[0~~\id_P]$. We claim that this is indeed an almost split sequence in $\mon(\omega, \cp)$. First, we show that $\eta$ is non-split. Assume on the contrary that $\eta$ is a split sequence. Hence, there is a morphism $\psi=(\psi_1, \psi_2):(P\st{f}\rt Q)\lrt(Q\oplus P\st{q}\rt P\oplus Q)$ such that $g\psi=\id_{(P\st{f}\rt Q)}$.  Since by Theorem \ref{ds}, the functor $T:\HT\mon(\omega, \cp)\lrt\ugp(R)$ is fully faithful, one may find a morphism $h:M\rt X$ in $\ugp(R)$ such that $T(\psi)=h$.  Consequently, $f_2h=T(g)T(\psi)=T(g\psi)=\id_{T(P\rt Q)}=\id_M$ in $\ugp(R)$, and so, $f_2h-\id_M$ factors through a projective $R$-module. Now Lemma \ref{split} forces $f_2$ to be a split epimorphism in $\gp(R)$, which is a contradiction, and then, $\eta$ is non-split.

Next assume that $\gamma=(\gamma_1, \gamma_0):(P'\st{f'}\rt Q')\lrt(P\st{f}\rt Q)$ is a non-split epimorphism in $\mon(\omega, \cp)$. We have to show that $\gamma$ factors through $g$. By our hypothesis, there is a commutative diagram of $S$-modules with exact rows; $$\xymatrix{0\ar[r] & P'\ar[r]^{f'} \ar[d]_{\gamma_1} &\ Q' \ar[d]_{\gamma_0}\ar[r] & \cok f' \ar[d]_{e}\ar[r]& 0 \\0\ar[r]&P \ar[r]^{f} &  Q \ar[r] & M\ar[r]& 0.}$$ A similar argument has been appeared above, reveals that $e$ is a non-split epimorphism. Now since  $0\rt\tau(M)\st{f_1}\rt X\st{f_2}\rt M\rt 0$ is an almost split sequence, we will have an $R$-homomorphism $k_1:\cok f'\rt X$ in which $f_2k_1=e$. Hence, another use of the fact that the functor $T$ is full, gives us a morphism $\alpha=(\alpha_1, \alpha_0):(P'\st{f'}\rt Q')\lrt(Q\oplus P\rt P\oplus Q)$ with $T(\alpha)=k_1$. Now the faithfulness of $T$ enables us to deduce that $g\alpha=\gamma$ in $\HT\mon(\omega, \cp)$, and so by Proposition \ref{rrrr}, $g\alpha-\gamma$ factors through a projective object of $\mon(\omega, \cp)$. Therefore, applying Lemma \ref{split}, gives us a morphism $\psi=(\psi_1, \psi_0):(P'\st{f'}\rt Q')\lrt(Q\oplus P\rt P\oplus Q)$ in which $g\psi=\gamma$. {Consequently, $g$ is a right almost split morphism. Moreover, thanks to the existence of an almost split sequence ending at $M$, one infers that $M$ is indecomposable, and so, applying
Proposition \ref{p1} yields that $(Q\st{f_{\si}}\rt P)$ has local endomorphism ring}. Thus, $\eta$ is an almost split sequence.\\
Case 2: $\di R$ is even. {So, in view of  Proposition \ref{dim} $\tau(M)= M$. Now one may follow the argument that appeared in  Case 1 and conclude that there is an almost split sequence in $\mon(\omega, \cp)$, ending at $(P\st{f}\rt Q)$, and in particular, $\tau_{\mon}(P\st{f}\rt Q)=(P\st{f}\rt Q)$.\\
Conversely, assume that there is an almost split sequence in $\mon(\omega, \cp)$ terminating in $(P\st{f}\rt Q)$. So, the above method gives rise to the existence of an almost split sequence ending at $M$.
 Thus the proof is completed.}
\end{proof}

Recall that  $R$ is called an isolated singularity, if $R_{\p}$ is a regular local ring, for all non-maximal prime ideals $\p$ of $R$. Moreover, an $R$-module $M$ is said to be locally projective on the punctured spectrum of $R$, provided that $M_{\p}$ is a projective $R_{\p}$-module, for all non-maximal prime ideals $\p$ of $R$.

Assume that $M$ is an indecomposable non-projective object of $\gp(R)$ which is locally projective on the punctured spectrum of $R$. So, in view of the main result of \cite{AR5} (see also \cite[Theorem 13.8]{leuschke2012cohen}), there is an almost split sequence ending at $M$ in $\gp(R)$. Particularly, if $R$ is an isolated singularity, then any non-projective indecomposable Gorenstein projective $R$-module appears as the right term of an almost split sequence in $\gp(R)$; see also \cite[Corollary 13.9]{leuschke2012cohen}. This fact together with Theorem \ref{ass}, enables us to quote the following results.

\begin{cor}Let $(P\st{f}\rt Q)$ be an object of $\mon(\omega, \cp)$ such that $\cok f$ is an indecomposable  non-projective object of $\gp(R)$ which is locally projective on the punctured spectrum of $R$. Then  $\tau_{\mon}((P\st{f}\rt Q))=(P\st{f}\rt Q)$, if $\di R$ is even, and $\tau_{\mon}((P\st{f}\rt Q))=(Q\st{f_{\si}}\rt P)$, otherwise. In particular, $\tau_{\mon}^2(f)=f$.
\end{cor}

\begin{cor}\label{cc}Let $R$ be  an isolated singularity. Then each indecomposable non-projective object of $\mon(\omega, \cp)$ appears as the right term of an almost split sequence.
\end{cor}
\begin{proof}Assume that $(P\st{f}\rt Q)$ is a non-projective indecomposable  object of $\mon(\omega, \cp)$. So by Propositon \ref{p1}(2), $\cok f$ is an indecoposable object of $\gp(R)$. {Moreover, as $(P\st{f}\rt Q)$ is non-projective, it is easily verified that the same is true for $\cok f$.} Now, since $R$ is an isolated singularity, as already declared above, by the main result of \cite{AR5} (also \cite[Corollary 13.9]{leuschke2012cohen}), there is an almost split sequence ending at $M$. Consequently, in view of Theorem \ref{ass}, there is an almost split sequence ending at $(P\st{f}\rt Q)$. So we are done.
\end{proof}

\begin{example}Assume that $\di S=1$. So $\di R=0$, meaning that $R$ is artinian. Since $R$ is self-injective, the categories $\md R$ and $\gp(R)$ are the same. In this case, it is known that  $\gp(R)$ admits almost split sequences; see \cite{auslander1975representation}. Consequently, by Corollary \ref{cc} each indecomposable non-projective object of $\mon(\omega, \cp)$ appears as the right term of an almost split sequence. Particularly, for a given non-projective indecomposable object $(P\st{f}\rt Q)\in\mon(\omega, \cp)$, we have $\tau_{\mon}(f)=f$.
\end{example}

%{\bf{Statements and Declarations.}} The authors declare that they have no known competing financial interests or personal relationships that could have appeared to influence the work reported in this paper.

%{\bf{Acknowledgments.}} The authors are grateful to the referee for reading the paper very carefully and giving a lot of valuable suggestions kindly and patiently.

\bibliographystyle{siam}
%\bibliography{newRefs(1)}

\end{document}